\documentclass[11pt,psfig]{article}
\usepackage{amssymb,amsmath,amscd}
\usepackage[all]{xy}
\newtheorem{lem}{Lemma}[section]
\newtheorem{thm}{Theorem}[section]
\newtheorem{cor}{Corollary}[section]

\newcommand{\f}[1]{\mathfrak{#1}}

\newcommand{\commentout}[1]{}

\newcommand{\mc}{\mathcal}
\newcommand{\arr}[1]{\left( \begin{array}{clcr} #1 \end{array} \right)}

\newcommand{\vertiii}[1]{{\left\vert\kern-0.25ex\left\vert\kern-0.25ex\left\vert #1 
    \right\vert\kern-0.25ex\right\vert\kern-0.25ex\right\vert}}

\begin{document}
\title{Representations of $ax+b$ group and Dirichlet Series}
\author{Hongyu He  \\
Department of Mathematics \\
Louisiana State University \\
email: hhe@lsu.edu\\
}
\date{}
\maketitle

\abstract{Let $G$ be the $ax+b$ group. There are essentially two irreducible infinite dimensional unitary representations of $G$, $(\mu, L^2(\mathbb R^+))$ and $(\mu^*, L^2(\mathbb R^+))$. In this paper, we give various characterizations about smooth vectors of $\mu$ and their Mellin transforms.   Let $\f d$ be a linear sum of delta distributions supported on the the positive integers $\mathbb Z^+$.   We study the Mellin transform of the matrix coefficients $\mu_{ \f d, f}(a)$ with $f$ smooth. We express these Mellin transforms  in terms of the Dirichlet series $L(s, \f d)$. We determine a sufficient condition  such that the generalized matrix coefficient $\mu_{\f d, f}$ is a locally integrable function and estimate the $L^2$-norms of $\mu_{\f d, f}$  over the Siegel set. We further derive an inequality which may potentially be used to study the Dirichlet series $L(s, \f d)$.   }

\section{Introduction}
Let $G$ be the $ax+b$ group, the semidirect product $\mathbb R \rtimes \mathbb R^+ $. Here $\mathbb R^+$ is the multiplicative group and $\mathbb R$ is the additive group. The group $G$ can be realized as the matrix group consisting of 
$$\arr{a & 0 \\ t  & 1}=\arr{ a & 0 \\ 0 & 1} \arr{1 & 0 \\ t & 1}, \qquad (a \in \mathbb R^+, t \in \mathbb R).$$
For simplicity, we denote the first factor by $a$ and second factor by $b_t$.
There are two equivalence classes of irreducible infinite dimensional unitary representations, $\mu$ and its contragredient $\mu^*$. The representation $\mu$ can be modeled on $L^2(\mathbb R^+)$ by defining
$$\mu(a b_t)f(x)=a^{-\frac{1}{2}} \exp (-2 \pi i t a^{-1} x) f(a^{-1} x).$$
The unitary representation $(\mu, L^2(\mathbb R^+))$  is well-known (\cite{fo}). Let $\mc H^{\infty}$ be the Frechet space of smooth vectors in $\mu$. The first result we proved in this paper is a characterization of smooth vectors in terms of their Mellin transform $\mc M$.
\begin{thm}
Let $g(z)$ be an analytic function on the right half plane $\mathbf H=\{ \Re(z) > \frac{1}{2} \}$. Then $g (z) \in \mc M (\mc H^{\infty})$ if and only if the following hold
\begin{enumerate}
\item For  any $n >0$ and  any finite closed interval $I$ in $( \frac{1}{2}, \infty)$, $g(\sigma+i s_1)$ is bounded by 
$C_{I, n} |\sigma+ i s_1|^{-n}$ for all $s_1 \in \mathbb R$ and $\sigma \in I$;
\item $g(\sigma+ i s_1)$ has a $L^2$ limit $g(\frac{1}{2}+i s_1)$ as $\sigma \rightarrow \frac{1}{2}$.
\item  For any $n \geq 1$, $ {s_1}^n g(\sigma+ i s_1)$ has a $L^2$ limit as $\sigma \rightarrow \frac{1}{2}$. 
\end{enumerate}
\end{thm}

Let $(\mc H^*)^{-\infty}$ be the dual topological vector space of $\mathcal H^{\infty}$, equipped with the weak star topology. We can extend  the representation $\mu$ to its dual topological vector space $(\mc H^*)^{-\infty}$. We obtain a representation $\mu$ on the distribution vectors.  Since  the Schwartz testing function $\mathbf D(\mathbb R^+)$ are dense in $\mc H^{\infty}$, (Theorem \ref{hinfinity}), the space $(\mc H^*)^{{-\infty}}$ is a linear subspace of $\mathbf D^{\prime}(\mathbb R)$ (\cite{tr}). Hence, all our discussion can be carried out in the classical sense.  \\
\\
Let $\f d= \sum_{n=1}^{\infty} d_n \delta_n$ be a distribution in $\mathbf D(\mathbb R^+)^{\prime}$. Here $\delta_n$ is the Dirac $\delta$ function supported on $n$ with $n$ a positive integer. Suppose that $\{ d_n \}$ grow slower than a polynomial. Then $\f d \in \mc H^{-\infty}$. In \cite{he}, we proved that
the classical definition of  matrix coefficients can be extended to a pair of distributions and the resulting generalized matrix coefficients are distributions on the group $G$. Hence generalized matrix coefficient $\mu_{\f d, f}(a b_t)$ is defined for any $f \in (\mc H^*)^{-\infty}$. 
\begin{thm}
 Let $f$ be a locally integrable function in $\mc H^{-\infty}$ such that for some $K >0$ and $\epsilon >0$  
$$\int_{K}^{\infty} a^{\epsilon} |f(a)|^2 d a < \infty.$$
Then $\mu_{\f d, f}(a b_t)$ is a well-defined locally integrable function. In addition, for some positive constant $C_{\epsilon, \f d}$
$$|d_1|^2 \int_{K}^{\infty} a^{\epsilon} |f(a)|^2 d a \leq \int_{K}^{\infty} \int_0^{1} |\mu_{\f d, f}(a b_t)|^2 a^{\epsilon} \frac{d a}{a} \leq C_{\epsilon, \f d} \int_{K}^{\infty} a^{\epsilon} |f(a)|^2 d a.$$
\end{thm}
Our estimate is the $L^2$-norm of $\mu_{\f d, f}(a b_t)$ over the Siegel set. \\
\\ 
Since $\{ d_n \}$ is of polynomial growth,  the Dirichlet series $L(s, \f d)= \sum_{n=1}^{\infty} d_n n^{-s}$ is well-defined  in a right half plane. For 
$\phi \in \mc (H^*)^{\infty}$, the matrix coefficient $\mu_{\f d, \phi}(b_T a)$ is a smooth function. It can be expressed in terms of Mellin transform.
\begin{thm}\label{mellinmatrix}
Suppose that the sequence $\{ d_n \}$ is bounded by a multiple of $n^{\epsilon}$ for any $\epsilon > 0$. Suppose that $L(s, \f d)$ is holomorphic and bounded by a polynomial on the strip
$\{ r \leq \Re(s) \leq 1+\epsilon \}$ for some $\epsilon > 0$. Then for any $\phi \in \mc H^{\infty}$, we have
$$\mu_{\f d, \phi}(b_T a)=\langle \mu(b_T a)\f d, \phi \rangle=\frac{1}{2 \pi} \int L(r+s_1 i, \f d) \mc M\{\phi(x) \exp - 2 \pi i Tx \} (r+s_1i) a^{-\frac{1}{2}-r-s_1 i} d s_1.$$
\end{thm}
The main result of this paper is the following inequality.
\begin{thm}
 Let $\f d=\sum d_n \delta_n$ with $\{ d_n \}$ bounded by $C_{\tau} n^{\tau}$ for any $ \tau >0$. Let  $T_1 \geq 1$, $\epsilon >0$ and $C_{\epsilon}=\sum_{n=1}^{\infty} n^{-1-\epsilon} |d_n|^2$. Let $f \in L^2(\mathbb R^+, (1+x^{\epsilon}) d x)$. We have
$$\frac{1}{2} T_1 C_{\epsilon} \int_1^{\infty} a^{\epsilon} |f(a)|^2 d a \leq \int_{1}^{\infty} \int_0^{T_1} |\mu_{\f d, f}(b_T a)|^2 d T a^{\epsilon} \frac{d a}{a} \leq 2 T_1 C_{\epsilon} \int_1^{\infty} a^{\epsilon} |f(a)|^2 d a.$$
Let $\chi_{[1, \infty)}$ be the indicator function of $[1 , \infty)$. Then there exists a constant $C_{\epsilon, T_1}$ such that
$$ \int_0^{T_1} \|\mc M^{\frac{\epsilon}{2}} \{ \mu_{\f d, f}( b_T a) \chi_{[1, \infty)}(a)\}(s_1)\|_{L^2(\mathbb R)}^2 d T  \leq  C_{\epsilon, T_1} \| \mc M^{\frac{1+\epsilon}{2}}(f \chi_{[, \infty)})(s_1) \|_{L^2(\mathbb R)}^2.$$
\end{thm}
Here $\mc M^{\frac{\epsilon}{2}}\{ * \}(s_1)$ is Fourier-Mellin transform, a variation of $\mc M \{ * \} (\frac{\epsilon}{2}+s_1 i)$.\\
\\
Notice that Mellin transform  $\mc M^{\frac{\epsilon}{2}} \{ \mu_{\f d, f}( b_T a) \chi_{[1, \infty)}(a) \} (s_1)$ may be obtained from $\mc M^{\frac{\epsilon}{2}} \{ \mu_{\f d, f}( b_T a)  \} (s_1)$ by solving a Riemann-Hilbert problem. By Theorem \ref{mellinmatrix},  our inequality may then be expressed in terms of $L(s, \f d)$. Furthermore, if we replace the $ax+b$ group by $SL(2, \mathbb R)$ and assume that $\{ d_n \}$ are the Fourier coefficients of an automorphic form, then the factor $\chi_{[1, \infty)}$ can be removed. We will be able to bound certain $L^2$-norm of  the Mellin transform
$$\mc M^{\frac{\epsilon}{2}} \{ \mu_{\f d, f}( g a)  \} (s_1)$$
by certain $L^2$-norm of the form $M^{\frac{1+\epsilon}{2}} \{f \}(s_1)$
(\cite{he1}). Tis will shed some lights on the behavior of the Dirichlet series $L(s, \f d)$.
 
\section{Representations of $ax+b$ group}
Let us define $ax+b$ group to be the semidirect product of $A \cong \mathbb R^+$ and $B \cong \mathbb R$. More precisely, for $a \in \mathbb R^+$ and $b_t \in \mathbb R$.
define $(a, b_t)=(a,0)(1, b_t)=(1, b_{a^{-1}t})(a,0)$. We may simply write $ab_t$ for $(a, b_t)$. Then
$$ab_t=b_{a^{-1}t} a.$$
The product rule is given by
$$ab_t a^{\prime} b_{t^{\prime}}= a a^{\prime} b_{a^{\prime} t+t^{\prime}}.$$
Let $G$ be the $ax+b$ group. The group $G$ can be parametrized by $a b_t$ or $b_{T} a$ with $T= a^{-1} t$.
The standard left invariant Haar measure is given by $ \frac{d a}{a} d t= d a d T$. The standard right invariant measure is given by
$\frac{d a}{a} d T= a^{-1} \frac{d a}{ a} dt$. \\
\\
The unitary dual of $G$ is easy to describe. There are one dimensional representations parametrized by  $i \lambda \in i \mathbb R$:
$$\chi_{i \lambda}: a b_t \rightarrow a^{i \lambda}.$$
There are two infinitely dimensional representations $(\mu, L^2(\mathbb R^+))$ and $(\mu^*, L^2(\mathbb R^+))$ defined as follows
$$\mu(a) f(x)=a^{-\frac{1}{2}} f(a^{-1} x), \qquad \mu^*(a) f(x)=a^{-\frac{1}{2}} f(a^{-1} x);$$
$$\mu(b_t)f(x)=\exp (-2 \pi i t x) f(x), \qquad \mu^*(b_t) f(x)=\exp (2 \pi i x t) f(x);$$
for  $f \in L^2(\mathbb R^+)$, where $ \mathbb R^+$ is equipped with the Euclidean measure. See for example \cite{fo} for details of the unitary dual of $G$. Notice that $\mu^*$ is simply the contragredient representation of $\mu$. The representations $\mu$ and $\mu^*$ differ by a complex conjugation on the action of $b_t$. We will mainly focus on $\mu$. The support of $\mu|_{B}$ is the negative half line $\mathbb R^-$ and the support of $\mu^*|_B$ is the positive half line $\mathbb R^+$. \\
\\
We start with the following lemma about the action of the Lie algebra $\f g$. 
\begin{lem}
The image of the universal enveloping algebra $\mu(U(\mathfrak g))$ is spanned by
$$\{ x^m \frac{d^n}{ d x^n}, m \geq n \geq 0 \}.$$
\end{lem}
Proof: Let $H$ be the infinitesimal generator of $A$ and $E$ be the infinitesimal generator of $B$. Then
$\mu(H)=-\frac{1}{2}- x \frac{d }{ d x}$ and $\mu(E)= - 2 \pi i x$.  Since the identity $1$ is contained in the universal enveloping algebra $U(\f g)$, $\mu(U(\f g))$ is generated by $\{x, x \frac{d}{d x} \}$ as a ring. One can then proceed by induction on $n$ to show that
$x^n \frac{d^n}{ d x^n} \in \mu(U(\mathfrak g))$ for any $n \geq 0$. It follows that 
$$\{ x^m \frac{d^n}{ d x^n}, m \geq n \geq 0 \} \subseteq \mu(U(\f g)).$$
Conversely, $\mu(U(\f g))$ is spanned by $x^k (x \frac{ d }{d x})^n$ since the Lie bracket $[x \frac{d }{d x}, x]=x$. But
$(x \frac{d}{ d x})^n$ is a linear combination of $x^k \frac{d^k}{ d x^k}$ with $0 \leq  k \leq n$. Hence
$$ \mu(U(\f g))   \subseteq \{ x^m \frac{d^n}{ d x^n}, m \geq n \geq 0 \}.$$
 $\Box$.

\subsection{smooth vectors}
Let $(\mu, \mathcal H=L^2(\mathbb R^+))$ be the irreducible unitary representation of $G$ defined above. Let $\mc H^{\infty}$ be the space of smooth vectors equipped with the canonical Frechet topology defined by seminorms $\{ \| D f \|; D \in \mu(U(\f g)) \}$. Sometimes the seminorm $\| D f \|$ is denoted by $\| f \|_D$.
\begin{lem} $$ \mc H^{\infty}=\{ f \in C^{\infty}(\mathbb R^+) \mid x^m \frac{d^n}{d x^n} f \in L^2(\mathbb R^+) \, \,\, \forall \, \, m \geq n \geq 0.\}.$$
\end{lem}
Proof: Let $f$ be a smooth vector. Then $D f \in L^2(\mathbb R)$ for any $D \in \mu(U(\mathfrak g))$. The higher order derivatives of $f$ are all locally $L^2$, hence locally $L^1$. Therefore any order derivatives of $f$ are continuous. This show that $f \in C^{\infty}(\mathbb R^+)$. $\Box$.\\
\\
We shall now give a more detailed description of $\mc H^{\infty}$.
\begin{thm}\label{hinfty}
$f(x) \in \mc H^{\infty}$ if and only if $ f \in C^{\infty}(\mathbb R^+)$ and the following are satisfied
\begin{enumerate}
\item For any $n \geq 0$, we have $ x^n \frac{d^n f}{d x^n} \in L^2(0,\delta)$ for some $\delta>0$;
\item For any $n, m \in \mathbb N$ and some $k>0$, there exists a constant $c_{m,n,k}$ such that 
$$ |\frac{d^n f(x)}{d x^n}|< c_{m,n,k} (1+|x|)^{-m} \qquad (x >k).$$
\end{enumerate}
\end{thm}
Roughly, the derivatives of $f$ must not behave too badly near zero and must be rapidly decaying near infinity. This follows from Sobelev inequality.  \\
\\
Proof: To prove "if" part, assume $(1)$ and $(2)$ and $m \geq n \geq 0$. From $(1)$, we have $ x^m \frac{d^n f}{d x^n} \in L^2(0,\delta)$. By the rapidly decaying condition $(2)$, $x^m \frac{d^n f}{d x^n} \in L^2(k,\infty)$. Finally, since $f \in C^{\infty}(\mathbb R)$, $x^m \frac{d^n f}{d x^n} $ will be bounded between $\delta$ and $k$, hence in $L^2(\delta,k)$. Therefore, $x^m \frac{d^n f}{d x^n} \in L^2(\mathbb R^+)$. \\
\\
Conversely, let $m \geq n \geq 0$, $(1)$ is obvious. Observe that
\begin{equation}
\begin{split}
& |a^m \frac{d^n f(a)}{d x^n}-\frac{d^n f(1)}{d x^n}|= |\int_1^a (x^m \frac{d^n f}{d x^n})^{\prime} d x| \\
=& |\int_1^a m x^{m-1} \frac{d^n f}{d x^n} + x^m \frac{d^{n+1} f}{d x^{n+1}} d x | \\
\leq & m \int_1^a | x^{m-1} \frac{d^n f}{d x^n} | d x+ \int_{1}^a |x^m \frac{d^{n+1} f}{d x^{n+1}}| d x \\
\leq & m (\int_1^a |x^{m} \frac{d^n f}{d x^n}|^2 dx)^{\frac{1}{2}} (\int_1^a |x^{-1}|^2 dx)^{\frac{1}{2}}+ (\int_1^a |x^{m+1} \frac{d^{n+1} f}{d x^{n+1}}|^2 d x )^{\frac{1}{2}}(\int_1^a |x^{-1}|^2 d x)^{\frac{1}{2}} \\
\leq & C_{n,m} \\
\end{split}
\end{equation}
Hence $x^m \frac{d^n f}{d x^n}(x)$ is bounded by a constant. $(2)$ follows. $\Box$ \\
\\
We shall observe that $C_c^{\infty}(\mathbb R^+) \subseteq \mc H^{\infty}$.

\subsection{2nd Characterization of $\mc H^{\infty}$}
For any $x \in \mathbb R^+$, let  $h= \ln x$. Define a unitary operator $\mc I: L^2(\mathbb R^{+}) \rightarrow L^2(\mathbb R)$ by
$$\mc I(f)(h)=(\exp \frac{1}{2} h) f(\exp h).$$
It is easy to check that $\mc I(\mu(H) f)(h)=-\frac{d}{d h}(\mc I(f))$ and $\mc I(\mu(E) f)(h)=(- 2 \pi i \exp h) \mc I(f)(h)$. Immediately, we have
\begin{lem}\label{phibound1}
$\phi \in \mc I(\mc H^{\infty})$ if and only if
\begin{enumerate}
\item $\phi \in C^{\infty}(\mathbb R)$;
\item for any $n \geq 0$, $\frac{ d^n \phi(h)}{d h^n} \in L^2(-\infty, k)$ for some $k$;
\item  $\frac{ d^n \phi(h)}{d h^n} $ is bounded by $C_{n,m}(1+\exp h)^{-m}$ for any $m>0$ and $h > 0$.
\end{enumerate}
\end{lem}
We need to analyze $(2)$ a little deeper. Since ${\phi}^{\prime} \in L^2(\mathbb R)$, for $h<-1$,
$$|\phi(-1)-\phi(h)|=|\int_{h}^{-1}{\phi}^{\prime}(y) d y| \leq  \sqrt{-1-h} \| {\phi}^{\prime}\|_{L^2}.$$
Hence $\phi(h)$ is bounded by a multiple of $\sqrt{|h|}$ at $-\infty$. If we work a little harder, we have
\begin{thm}\label{phibound2} If $\phi(h) \in \mc I(\mc H^{\infty})$, then
$ \phi(h)$ and all its higher order derivatives are  bounded.
\end{thm}
Proof: Consider the Fourier transform
$$\mc F\phi(\xi)=\int \phi(h) \exp (-2\pi i h \xi) d h.$$
Since $\phi^{(n)}(h) \in L^2(\mathbb R)$, we have $\xi^n \mc F \phi(\xi) \in L^2(\mathbb R)$ for any $n \geq 0$. Then
$$|\phi(h)|=|\int \mc F \phi(\xi) \exp (2 \pi i h \xi) d \xi| \leq \int_{-\infty}^{-1} |\mc F \phi(\xi)| d \xi+\int_{1}^{\infty} |\mc F \phi(\xi)|d \xi +C.$$
The first two terms can all be bounded by the $L^2$-norm of $\xi \mc F \phi(\xi)$. Hence $\phi(h)$ is bounded. Similarly,
$\frac{d^n \phi}{d h^n}(h)$ is bounded for all $n \geq 0$. $\Box$ \\
\\
This estimate applies to all $\phi(h)=\exp \frac{1}{2} h \, f(\exp h)$. Returning to $f$, we have
\begin{cor}
Let $f \in \mc H^{\infty}$. Then all the derivatives of $ f$ will be bounded by a multiple of $x^{-\frac{1}{2}}$ near $0$.
\end{cor}

\section{Fourier-Mellin transform of $\mc H^{\infty}$}
Let $f$ be a locally integrable function on $\mathbb R^+$. Narrowly speaking, Mellin transform $$\mc M f(s)= \int f(x) x^s \frac{d x}{x}$$
 is required to have a band of convergence. In this paper, we will loosely define
$\mc M f(s)=\int f(x) x^s \frac{d x}{x}$ as long as the integral converges absolutely. Indeed, if $f(x) x^{\sigma-1} \in L^1(\mathbb R^+)$ we can define
$$\mc M^{\sigma} f(s_1)=\int f(x) x^{\sigma+i s_1} \frac{d x}{x}.$$
Throughout this paper $s_1$ will always be real and $s $ will always be complex with its imaginary part $i s_1$.\\
\\
If $x=\exp h$, then $$\mc M^{\sigma} f(s_1)=\int f( \exp h) \exp (\sigma+i s_1) h d h.$$
  Notice that $\mc M^{\sigma}$ can be extended to $L^2$-functions and beyond, by Fourier transform. We call $\mc M^{\sigma} $ the Fourier-Mellin transform. In this context, the Fourier-Mellin inversion for $L^2$-functions will be given by
$$f(x)=\frac{1}{2 \pi} \int M^{\sigma}f(s) x^{-\sigma- s_1 i} d s_1. $$
 Unless such interpretation is needed, we will generally retain the usage of Mellin-transform  $\mc M f(s)$ or Mellin inversion.

\begin{thm}\label{minfty} Suppose that $f \in \mc H^{\infty}$. Then 
\begin{enumerate}
\item $\mc M f(s)$ is defined and holomorphic on the open half plane $\mathbf H=\{ \Re s > \frac{1}{2}\}$;
\item $\mc M f(s)$ decays faster than any $|s|^{-n}$ along any compact vertical stripe in  $\mathbf H$;
\item $\mc M^{\sigma} f(s_1)$ has an $L^2$-limit as $\sigma \rightarrow \frac{1}{2}$. The $L^2$-limit is the Fourier-Mellin transform $\mc M^{\frac{1}{2}} f(s_1)$.
\end{enumerate}
\end{thm}
Proof: Since $f$ is bounded by a multiple of $x^{-\frac{1}{2}}$ near $0$ and $x^{-n}$ near $\infty$, $\mc M f(s)$ is well-defined and analytic on the open half plane $\mathbf H$. $(1)$ is proved. \\
\\
Observe that 
$$\mc M f(\sigma+ i s_1)=\int \mc I (f)( h) \, \exp (\sigma-\frac{1}{2})h \, \exp is_1 h \, d h$$ in terms of Fourier transform. By Lemma \ref{phibound1} and Theorem \ref{phibound2}, if $\sigma> \frac{1}{2}$,
$\mc I(f)(h) \exp (\sigma-\frac{1}{2})h$ and all its derivatives will be in $L^1$. It follows that $\mc M f(\sigma+ i s_1)$ decays faster than any $|s_1|^{-n}$ as $|s_1| \rightarrow \infty$ for $\sigma$ in a finite closed interval in $(\frac{1}{2}, \infty)$. $(2)$ is proved. \\
\\
Let $\sigma > \frac{1}{2}$. Observe that $\mc I(f)(h) \exp (\sigma-\frac{1}{2})h \rightarrow \mc I(f)(h)$ in $L^2$ as $\sigma \rightarrow \frac{1}{2}$. On the Fourier-Mellin transform side, we have $\mc M^{\sigma} f(s_1) \rightarrow \mc M^{\frac{1}{2}} f (s_1)$ in $L^2(\mathbb R)$. $\Box$.\\
\\
We summarize how the Lie algebra act on the Mellin transform $\mc M f(s)$.
\begin{lem}\label{mliealgebra} Let $f \in \mc H^{\infty}$. Then
$\mc M (xf)(s)= \mc M (f)(s+1)$ and $\mc M(\mu(H) f)(s)=(s-\frac{1}{2}) \mc M(f)(s)$ for any $s \in \mathbf H$.
\end{lem}
Notice that Theorem \ref{minfty} (1)(2) are preserved under the action of Lie algebra. We can now characterize the Mellin transform of $\mc H^{\infty}$.
\begin{thm}
Let $g(s)$ be an analytic function on $\mathbf H$. Then $g (s) \in \mc M (\mc H^{\infty})$ if and only if the following hold
\begin{enumerate}
\item For  any $n >0$ and  any finite closed interval $I$ in $( \frac{1}{2}, \infty)$, $g(\sigma+i s_1)$ is bounded by 
$C_{I, n} |\sigma+ i s_1|^{-n}$ for all $s_1 \in \mathbb R$ and $\sigma \in I$;
\item $g(\sigma+ i s_1)$ has a $L^2$ limit $g(\frac{1}{2}+i s_1)$ as $\sigma \rightarrow \frac{1}{2}$.
\item  For any $n \geq 1$, $ {s_1}^n g(\sigma+ i s_1)$ has a $L^2$ limit as $\sigma \rightarrow \frac{1}{2}$. 
\end{enumerate}
\end{thm}
Proof: The only if part can be easily seen from Theorem \ref{minfty} and Lemma \ref{mliealgebra}. To show "if" part, assume $(1)(2)(3)$. We must show that $g(z)= \mc M f(z)$ for some $f(x) \in \mc H^{\infty}$. First the $L^2$-limit of $ {s_1}^n g(\sigma+ i s_1)$ will be $ {s_1}^n g(\frac{1}{2}+ i s_1)$ as $\sigma \rightarrow \frac{1}{2}$ for $s_1$ in any compact interval. Hence $ {s_1}^n g(\sigma+ i s_1) \rightarrow  {s_1}^n g(\frac{1}{2}+ i s_1)$ in $L^2(\mathbb R)$ as $ \sigma \rightarrow \frac{1}{2}$. \\
\\
Let
$$ f(x)= \frac{1}{2 \pi} \int g(\sigma+s_1 i) x^{-\sigma-i s_1} d s_1 \qquad (\sigma>\frac{1}{2}).$$
By (1), the integral above is independent of $\sigma$. In addition, for any $n \in \mathbb N$, 
$$\frac{1}{2 \pi} \int g(\sigma+s_1 i) [\prod_{j=0}^{n-1} (-\sigma-i s_1-j) ] x^{-\sigma-i s_1-n} d s_1$$
 converges absolutely and defines the higher order derivative $\frac{d^n f}{d x^n}(x)$. Hence $f(x) \in C^{\infty}(\mathbb R^+)$. By letting $\sigma$ arbitrarily large, $\frac{d^n f}{d x^n}(x)$ is uniformly bounded by a multiple of $|x|^{-m}$ for $x>1$. This proves $(2)$ of Theorem \ref{hinfty}. Near $0$, by letting $\sigma \rightarrow \frac{1}{2}$, $f(x)$ is bounded by a multiple of $x^{-\frac{1}{2}-\delta}$ for arbitrarily small $\delta>0$. It follows that $\mc Mf(s)$ is well-defined on $\mathbb H$. By Mellin inversion, $\mc Mf(s)=g(s)$ for any $s \in \mathbb H$.\\
\\
To show that $f \in \mc H^{\infty}$,  it suffices to show $(1)$ of Theorem \ref{hinfty}, namely $x^n \frac{d^n f}{d x^n} \in L^2(0,1)$. By Fourier-Mellin inversion, $x^{\sigma} f(x) \in L^2(\mathbb R^+, \frac{d x}{x})$ and $|x^{\sigma-\frac{1}{2}}f(x)| \in L^2(\mathbb R^+)$ for any $\sigma > \frac{1}{2}$. By  $(2)$, 
$|x^{\sigma-\frac{1}{2}}f(x)|$ has a $L^2$ limit in $L^2(\mathbb R^+)$ as $\sigma \rightarrow \frac{1}{2}$. This $L^2$-limit must be the pointwise limit $f(x)$. Hence $f(x) \in L^2(\mathbb R^+)$. Similarly, by $(3)$ and Lemma \ref{mliealgebra}, $(x \frac{d }{d x}+\frac{1}{2})^n f \in L^2(\mathbb R^+)$. By Theorem \ref{hinfty}, $f \in \mc H^{\infty}$. $\Box$ \\
\\

\section{ Dirichlet series and distributions supported on $\mathbb Z^+$}
Recall that $C_c^{\infty}(\mathbb R^+) \subseteq \mc H^{\infty}$. In addition, for any compact interval $K \subset \mathbb R^+$, the seminorms of $\mc H^{\infty}$ yield the seminorms of $C_c^{\infty}(K)$ when restricted to $C_c^{\infty}(K)$. Thus every continuous functional on the Frechet space $\mc H^{\infty}$ is a distribution on $\mathbb R^+$ in the sense of Schwartz. However a distribution on $\mathbb R^+$ may not be a distribution on $\mc H^{\infty}$. Let $\mc H^{-\infty}$ be the space of continuous linear functionals on $\mc H^*$, equipped with the weak-star topology.\\
\\
Let $d=(d_1,d_2, \ldots, d_n \ldots)$ be an infinite sequence of complex numbers. We say that $d$ is of polynomial growth if there exists a $k \in \mathbb N$ and $C_k>0$ such that $d_n < C_k n^k$ for all $n \in \mathbb N$. 
Let $\delta_n$ be the $\delta$ distribution on $n \in \mathbb N$. Define the distribution $\mathfrak d= \sum_{n=1}^{\infty} d_n \delta_n$. If $d$ is of polynomial growth, define  the Dirichlet series $L(s, \mathfrak d)=\sum_{n=1}^{\infty} d_n n^{-s}$. $L(s, \mathfrak d)$ is well-defined for $\Re(s)> k+1$.
\subsection{Distributions supported on $\mathbb N$}
\begin{thm} If $d$ is of at most polynomial growth, then $\f d \in \mc H^{-\infty}$.
\end{thm}

Proof:  It suffices to show that there exists a $F \in \mc H$ and $D \in U(\f g)$ such that $\mathfrak d= \mu(D) F$. 
\commentout{
We will work with the model $\mc I(\mc H^*)$. Recall $\mc I: (\mc H^*)^{\infty} \rightarrow \mc I((\mc H^*)^{\infty}) \subseteq L^2(\mathbb R)$. The dual isomorphism is given by $\mc I^*: (\mc I((\mc H^*)^{\infty}))^* \rightarrow \mc H^{-\infty}$. We are interested in $(\mc I^*)^{-1}$ that maps distributions
in $\mc H^{-\infty}$ to continuous linear functionals on $\mc I((\mc H^*)^{\infty})$. The space $(\mc H^*)^{\infty}$ is the same as $\mc H^{\infty}$. The action of $E \in \mathfrak g$ differs only by a sign. For the purpose of the proof of our theorem, we shall ignore this difference. We will simply speak of
$(\mc I^*)^{-1}: \mc H^{-\infty} \rightarrow \mc I(\mc H^{\infty})^{*}$. 

\begin{lem}
$(\mc I^*)^{-1} \delta_n =\frac{1}{\sqrt{n}} \delta_{\ln n}$.
\end{lem}
 }
 Consider the antiderivative of $\f d$:
$$f(x)=\sum_{n \leq x} d_n \qquad (x >0).$$
Obviously, $f(x) \leq C_k |x|^{k+1}$ for some $C_k > 0$ and $k \in \mathbb N$. Let $F(x)=x^{-k-2} f(x)$. Then $F(x) \in \mc H$ and
$$\mathfrak d = \frac{d}{d x} x^{k+2} F(x)=x^{k+2} \frac{d}{d x} F(x)+(k+2) x^{k+1} F(x)=(x^{k+2} \frac{d}{d x}+(k+2) x^{k+1}) F(x).$$
Since $x^{k+2} \frac{d}{d x}+(k+2) x^{k+1} \in \mu(U(\f g))$, $\f d \in \mc H^{-\infty}$. $\Box$ \\
\\

The Dirichlet series $L(s, \f d)$ is well-defined for $s>k+1$ if $\{ d_n \}$ is bounded by $\{ C_k n^{k} \}$ uniformly. 
We say that $L(s, \f d) \in \mc R$ if $\{ d_n \}$ is bounded by $C_{\epsilon} n^{\epsilon}$ for any $\epsilon >0$. This is the Ramanujan condition. 
Clearly, $L(s, \f d)$ is an analytic function on the open half plan $\mathbb H_1= \{z \in \mathbb C \mid \Re(z) > 1 \}$. In addition, if $L(s, \f d) \in \mc R$ and $L(s, \f d^{\prime}) \in \mc R$, then $L(s, \f d) L(s, \f d^{\prime})$ is also analytic on $\mathbb H_1$. In fact $L(s, \f d) L(s, \f d^{\prime}) \in \mc R$. The following is well-known. 
\begin{lem} The product $L(s, \f d) L(s, \f d^{\prime})=L(s, \mathfrak j)$
with
$$j_l=\sum_{mn=l} d_m d_n^{\prime}.$$
If $L(s, \f d), L(s, \f d^{\prime}) \in \mc R$, then $L(s, \f d) L(s, \f d^{\prime}) \in \mc R$. 
\end{lem}

Hence $\mc R$ is a ring. We write $j$ as $d \sharp d^{\prime}$ and $\f j$ as $\f d \sharp \f d^{\prime}$.
\commentout{
Proof: Fix any $\delta=\frac{\epsilon}{2}$. Assume that $|d_n| \leq C_{\delta} n^{\delta}$ and $|d_m^{\prime}| \leq C_{\delta} m^{\delta}$. Then
$$j_k \leq C_{\delta}^2 k^{\delta} \#({\rm factors }\,\, of \,\, k).$$
We can prove that $\#({\rm factors} \,\, of \,\, k) \leq 2 \log_2 k$ as follows. Let $k=\prod_{i=1}^l p_i^{\alpha_i}$ be the prime factorization of $k$.
Then 
\begin{equation}
\begin{split}
\#(factors \,\, of \,\, k) = & \sum_{i=1}^l (\alpha_i+1) \\
 = & \sum_{i=1}^l 1+\log_{p_i} p_i^{\alpha_i} \\
 \leq & \sum_{i=1}^l 2 \log_2 p_i^{\alpha_i} \\
 = & 2 \log_2 k
\end{split}
\end{equation}
Then $\{ j_k \}$ will be bounded uniformly by $2 C_{\delta}^2 k^{\delta} \log_2 k < C_{\epsilon} k^{\epsilon}$ for any $k$. $\Box$
}
\subsection{$\langle \f d, f \rangle$ in terms of Mellin transform}
Let $\mc{MR}$ be the subring of $\mc R$ consisting of $g(z) \in \mc R$ with meromorphic continuation to $\mathbb C$. We shall see how the pairing
$\langle \f d, f \rangle$ can be expressed in terms of $L(s, \f d)$.
\begin{lem} Let $f \in \mc H^{\infty}$ and $L(s, \f d) \in \mc{MR}$. Then for any $ r >1$, we have
$$\langle \f d, f \rangle=\frac{1}{2 \pi} \int L(r+s_1i, \f d) \mc Mf(r+s_1i) d s_1.$$
\end{lem}
Proof: Let $r >1$. By the characterization of Mellin transform of $\mc H^{\infty}$ and Mellin inversion, we have
$$f(x)= \frac{1}{2 \pi} \int \mc Mf(r+ i s_1) x^{-r- i s_1} d s_1$$
with $\mc Mf(r+s_1 i)$ fast decaying. Due to absolute convergence of $L(r+s_1 i, \f d)=\sum_{k=1}^{\infty} \frac{d_k}{k^{r+s_1 i}}$, we have
$$\langle \f d, f \rangle= \sum_{k=1}^{\infty} d_k f(k)=\frac{1}{2 \pi} \sum_{k=1}^{\infty} d_k \int \mc Mf(r+ i s_1)  {k^{-r-s_1 i}} d s_1$$
$$=\frac{1}{2 \pi} \int \mc Mf(r+ i s_1) \sum_{k=1}^{\infty} d_k  k^{-r-s_1 i} d s_1 =\frac{1}{2 \pi} \int L(r+s_1 i, \f d) \mc Mf(r+s_1i) d s_1.$$
Now we can move the integral from $r+ i \mathbb R$ to the left of $1+  i\mathbb R$.
\begin{thm}\label{mellinpair}
Let $f \in \mc H^{\infty}$, $1 \geq r > \frac{1}{2}$ and $L(s, \f d) \in \mc{MR}$. Suppose that $L(s, \f d)$ is holomorphic and bounded by a polynomial on the strip
$\{ r \leq \Re(s) \leq 1+\epsilon \}$ for some $\epsilon > 0$. Then
$$\langle \f d, f \rangle=\frac{1}{2 \pi} \int L(r+s_1i, \f d) \mc Mf(r+s_1i) d s_1.$$
\end{thm}
We shall make some remarks here. First this result can be easily extended to the situation that $L(s, \f d)$ has finite number of poles in the stripe $\{ r \leq \Re(s) \leq 1+\epsilon \}$ by including the residue at the poles in the equation. Care should be taken to manage the distribution of infinitely many poles. But for the applications to $L$-functions, this theorem is sufficient. Secondly, this result can also be extended to the left half plane for $L$-functions. However, in our context, the Mellin transform $\mc M f$ is only defined on the right half plane. Hence there is a nature barrier at the critical line. \\
\\
For any $\sigma > \frac{1}{2}$, $x^{\sigma-1+i s_1}$ is a distribution in $\mc H^{-\infty}$. We can see this by writing $x^{\sigma-1+is_1}$ as
$(x+1) \frac{x^{\sigma+i s_1-1}}{x+1}$ and notice that $x+1$ is an operator in $\mu(U(\f g))$. 
\begin{cor} Let $ r > \frac{1}{2}$ and $L(s, \f d) \in \mc M \mc R$.  Suppose that $L(s, \f d)$ is holomorphic and bounded by a polynomial on the strip
$\{ r \leq \Re(s) \leq 1+\epsilon \}$ for some $\epsilon > 0$. Then as a distribution in $\mc H^{-\infty}$,
$$\f d = \frac{1}{2 \pi} \int^* L(r+ s_1 i, \f d) x^{r+s_1 i-1} d s_1.$$
The right hand side is a weak-star integral with test functions in $\mc H^{\infty}$.
\end{cor}
\subsection{$\langle \mu(a) \f d, f \rangle$ in terms of Mellin transform}
Observe that
$$\mc M (\mu^*(a^{-1})f)(s)=\int_0^{\infty} a^{\frac{1}{2}} f( a x) x^s \frac{d x}{x}=a^{\frac{1}{2}-s} \mc Mf(s).$$
\begin{cor} Under the same assumption as Theorem \ref{mellinpair}, 
$$\langle \mu(a) \f d, f \rangle= \frac{1}{2 \pi} \int L(r+s_1i, \f d) \mc Mf(r+s_1i) a^{\frac{1}{2}-r-s_1 i} d s_1.$$
or equivalently
$$\mc M (\langle \mu(a) \f d, f \rangle)(r-\frac{1}{2}+ s_1 i)= L(r+s_1 i, \f d) \mc M f(r+s_1 i).$$
\end{cor}

\section{Fundamental Inequality}
We shall now consider the (generalized) matrix coefficient $\langle \mu(a b_t) \f d, f \rangle$ with $ f \in (\mc H^*)^{\infty}$. We may still regard $f$ as a function in $\mc H^{\infty}$ and keep in mind that the action of $G$ will be given by $\mu^*$. We have
$$\langle \mu(a b_t) \f d, f \rangle= \langle \mu(b_t) \f d, \mu^*(a^{-1}) f \rangle=\sum_{n \in \mathbb Z^+} a^{\frac{1}{2}} d_n f(an) \exp (-2\pi i tn).$$
Since $f(x)$ decays very fast, $\sum_{n \in \mathbb N}$ converges as a Fourier series.
We write this matrix coefficient as $\mu_{\f d, f}(a b_t)$.
\subsection{Fourier Series Estimates}
By Parseval's theorem, we have
\begin{lem} Let $\f d$ be of polynomial growth and $f \in \mc H^{\infty}$, we have
$$\int_{0}^1 | \mu_{\f d, f}(a b_t)|^2 d t = a \sum_{n=1}^{\infty}  |f(a n) d_n |^2.$$
\end{lem}
\begin{thm}\label{fbd}
Let $L(s, \f d) \in \mc R$, $f \in \mc H^{\infty}$,  and $\epsilon >0$. Set  $C_{\epsilon}=\sum_{n=1}^{\infty} n^{-1-\epsilon} |d_n|^2$. Then for any $\delta \geq 0$, we have 
$$|d_1|^2 \int_{\delta}^{\infty} a^{\epsilon} |f(a)|^2 d a \leq \int_{\delta}^{\infty} \int_0^{1} |\mu_{\f d, f}(a b_t)|^2 a^{\epsilon} \frac{d a}{a} \leq C_{\epsilon} \int_{\delta}^{\infty} a^{\epsilon} |f(a)|^2 d a.$$

\end{thm}
Proof: Since $\f d$ satisfies the Ramanujan condition,  $C_{\epsilon}$ is well-defined. The lower bound follows from 
$$\int_0^1 | \mu_{\f d, f}( a b_t)^2| d t \geq a^{\epsilon} | f (a) d_1|^2.$$ 
The upper bound can be established as follows. 
\begin{equation}
\begin{split}
 & \int_{\delta}^{\infty} \int_0^{1} |\mu_{\f d, f}(a b_t)|^2 t a^{\epsilon} \frac{d a}{a} \\
 = & \int_{\delta}^{\infty} \sum_{n=1}^{\infty} |f(an) d_n |^2 a^{\epsilon} d a \\
 = & \sum_{n=1}^{\infty} (|d_n|^2 \int_{\delta}^{\infty} a^{\epsilon} |f(an)|^2  d a )\\
 = & \sum_{n=1}^{\infty} (|d_n|^2 \int_{n \delta}^{\infty} a^{\epsilon} n^{-1-\epsilon} |f(a)|^2 d a )\\
 \leq  & (\sum_{n=1}^{\infty} n^{-1-\epsilon} |d_n|^2)  \int_{{\delta}}^{\infty} a^{\epsilon} |f(a)|^2 d a.
 \end{split}
\end{equation}
 $\Box$.\\
\\
Notice that the domain we worked with $\{ \delta < a, 0 \leq t \leq 1 \}$ is classically known as the Siegel domain. However, the measure we use is different. 
In particular, if $\delta=0$, we have
$$|d_1|^2 \int_{\mathbb R^+} a^{\epsilon} |f(a)|^2 d a \leq \int_{\mathbb R^+} \int_0^{1} |\mu_{\f d, f}(a b_t)|^2 a^{\epsilon} \frac{d a}{a} \leq C_{\epsilon} \int_{\mathbb R^+} a^{\epsilon} |f(a)|^2 d a.$$

\subsection{Generalized Matrix Coefficients as functions}
Given any $f \in \mc H^{-\infty}$, we have shown in \cite{he} that generalized matrix coefficient $\langle \mu(g) \f d, f \rangle$ is a distribution on $G$.
The question then arises, for which $f$, $\langle \mu(g) \f d, f \rangle$ will be a locally integrable function. One simple guess is $f \in \mc H$? This turns out not to be the case for general $\f d$ satisfying the Ramanujan condition (\cite{zyg}). However, if we make assumption that for some $\epsilon >0$
and $K >0$,$$\int_{K}^{\infty} a^{\epsilon} |f(a)|^2 d a < \infty$$
 then $\mu_{\f d, f}(a b_t)$ will be a locally integrable function. 
We first prove a weaker version of this.
\begin{thm}\label{ine}
Let $f$ be a locally square integrable function in $\mc H^{-\infty}$ such that for some $\epsilon >0$ ,
$$\int_{0}^{\infty} a^{\epsilon} |f(a)|^2 d a < \infty.$$
Then $\mu_{\f d, f}(a b_t)$ is a well-defined locally square integrable function. In addition,
$$|d_1|^2 \int_{0}^{\infty} a^{\epsilon} |f(a)|^2 d a \leq \int_{0}^{\infty} \int_0^{1} |\mu_{\f d, f}(a b_t)|^2 a^{\epsilon} \frac{d a}{a} \leq C_{\epsilon} \int_{0}^{\infty} a^{\epsilon} |f(a)|^2 d a.$$
\end{thm}
Proof: Let $\{ f_n \}$ be a sequence of function in $C_{c}^{\infty}(\mathbb R^+) \subseteq \mc H^{\infty}$ such that
$$\int |f_n(a)-f(a)|^2 a^{\epsilon} d a \rightarrow 0 \qquad (as \, \, n \rightarrow \infty).$$
In particular, $\{ f_n \}$ is a Cauchy sequence in $L^2(\mathbb R^+, x^{\epsilon} d x)$.
By  Theorem \ref{fbd}, $\{ \mu_{\f d, f_n}(a b_t) \}$ is a Cauchy sequence in $L^2(G, a^{\epsilon} d t \frac{d a}{a})$.  It has a $L^2$-limit $\psi$  under the measure $a^{\epsilon} d t \frac{d a}{a}$.  Hence, as distributions in $\mathbf D(G)^{\prime}$,  $\mu_{\f d, f_n}(g) \rightarrow \psi(g)$ under the weak star topology.  By \cite{he}, $\mu_{\f d, f}(g)=\psi(g) \,\, a.e.$. For every $n$, 
$$|d_1|^2 \int_{0}^{\infty} a^{\epsilon} |f_n(a)|^2 d a \leq \int_{0}^{\infty} \int_0^{1} |\mu_{\f d, f_n}(a b_t)|^2 a^{\epsilon} \frac{d a}{a} \leq C_{\epsilon} \int_{0}^{\infty} a^{\epsilon} |f_n(a)|^2 d a.$$
By taking $n \rightarrow \infty$, our inequalities follow. $\Box$\\
\\
Knowing that $\mu_{\f d, f}(a b_t)$ is a square integrable function, by \cite{he}, we can approximate $\mu_{\f d, f}(a b_t)$ by $\mu_{\f d, f_n}( ab_t)$ as long as $f_n \rightarrow f$ under the weak star topology. 

\begin{thm}\label{ine1}
Let $f$ be a locally integrable function in $\mc H^{-\infty}$ such that for some $K >0$ and $\epsilon >0$  
$$\int_{K}^{\infty} a^{\epsilon} |f(a)|^2 d a < \infty.$$
Then $\mu_{\f d, f}(a b_t)$ is a well-defined locally integrable function. In addition,
$$|d_1|^2 \int_{K}^{\infty} a^{\epsilon} |f(a)|^2 d a \leq \int_{K}^{\infty} \int_0^{1} |\mu_{\f d, f}(a b_t)|^2 a^{\epsilon} \frac{d a}{a} \leq C_{\epsilon} \int_{K}^{\infty} a^{\epsilon} |f(a)|^2 d a.$$
\end{thm}
Proof: We may write $f$ as a sum of two locally integrable functions $f_1+f_2$ with $f_1 (x)=0 \,\, (\forall x \geq K)$ and $f_2(x)=0 \,\, (\forall x < K)$. Then $\mu_{\f d, f}(ab_t)= \mu_{\f d, f_1}(ab_t) + \mu_{\f d, f_2}(ab_t)$. By Theorem \ref{ine}, $ \mu_{\f d, f_2}(ab_t)$ is well-defined, locally square integrable. Since $f_1(x)=0$ for any $x \geq K$,  $\mu_{\f d, f_1}(a b_t)=\sum_{n \in \mathbb Z^+} a^{\frac{1}{2}} d_n f(an) \exp (-2\pi i tn)$ is locally a finite sum of locally integrable function, hence locally integrable. In addition, $\mu_{\f d, f_1}(a b_t)=0$ for any $a \geq K$. It follows that $ \mu_{\f d, f}(ab_t)=\mu_{\f d, f_2}(ab_t)$ for any $a \geq K$. By essentially the same proof as Theorem \ref{ine}, 
$$|d_1|^2 \int_{K}^{\infty} a^{\epsilon} |f(a)|^2 d a \leq \int_{K}^{\infty} \int_0^{1} |\mu_{\f d, f}(a b_t)|^2 a^{\epsilon} \frac{d a}{a} \leq C_{\epsilon} \int_{K}^{\infty} a^{\epsilon} |f(a)|^2 d a.$$
$\Box$
\subsection{Fundamental Inequality}
Define  the subset of $G$,  $\Omega_{T_1}= \{ b_T a \mid T \in [0, T_1], a \in [1, \infty) \}$. In terms of $ab_t$ coordinates, we have
$$\Omega_{T_1}= \{ a b_t \mid t \in [0, a T_1], a \in [1, \infty) \}.$$ 
\begin{thm}\label{fund} Let $L(s, \f d) \in \mc R$,  $T_1 \geq 1$ and $\epsilon >0$. Let $f \in L^2(\mathbb R^+, (1+x^{\epsilon}) d x)$. Set  $C_{\epsilon}=\sum_{n=1}^{\infty} n^{-1-\epsilon} |d_n|^2$. We have
$$\frac{1}{2} T_1 C_{\epsilon} \int_1^{\infty} a^{\epsilon} |f(a)|^2 d a \leq \int_{1}^{\infty} \int_0^{T_1} |\mu_{\f d, f}(b_T a)|^2 d T a^{\epsilon} \frac{d a}{a} \leq 2 T_1 C_{\epsilon} \int_1^{\infty} a^{\epsilon} |f(a)|^2 d a.$$
\end{thm}
Proof: If  $f(x) \in L^2(\mathbb R^+, (1+x^{\epsilon}) d x)$, then Theorem \ref{ine} holds.
To prove the upper bound, we observe
\begin{equation}
\begin{split}
 & \int_{1}^{\infty} \int_0^{T_1} |\mu_{\f d, f}(b_T a)|^2 d T a^{\epsilon} \frac{d a}{a} \\
 = & \int_{1}^{\infty} \int_{0}^{aT_1} |\mu_{\f d, f}(a b_t)|^2 a^{-1} d t a^{\epsilon} \frac{d a}{a} \\
 \leq & \int_1^{\infty} \lceil a T_1 \rceil a^{\epsilon-1} \int_{0}^1 |\mu_{\f d, f}(a b_t)|^2 d t \frac{d a}{a}\\
 \leq & \int_1^{\infty} 2 a T_1 a^{\epsilon-1} \int_{0}^1 |\mu_{\f d, f}(a b_t)|^2 d t \frac{d a}{a} \\
 = & 2 T_1 \int_{1}^{\infty} a^{\epsilon} \int_{0}^1 |\mu_{\f d, f}(a b_t)|^2 d t \frac{d a}{a} \\
 \leq & 2 T_1 C_{\epsilon} \int_{1}^{\infty} a^{\epsilon} |f(a)|^2 da 
\end{split}
\end{equation}
Here $\lceil a T_1 \rceil$ is the minimal integer greater or equal to $ aT_1$. We use the fact $\lceil a T_1 \rceil \leq 2 a T_1$.
The lower bound follows similarly. $\Box$\\
\\
If $T_1$ is any positive real number, the inequalities still hold if we replace $T_1 C_{\epsilon}$ be a constant $C_{T_1} C_{\epsilon}$. In addition, we can allow $a \in [\delta, \infty)$, then the inequalities will depend on $\delta$. However, we will not have any similar bounds if $a \in (0, \infty)$ since the behavior of $f(ab_t)$ near $a=0$ is hard to control without any additional assumption.  Nevertheless, If the group is $SL(2)$ and $\f d$ comes from the Fourier coefficients of an automorphic form, we can bound the $L^2$-norm of $\mu_{\f d, f}(B_T a)$ over the whole space $a \in (0, \infty)$ (\cite{he1}). 

\subsection{Inequality in terms of Mellin Transform}

For any locally integrable function $u(x)$ on $\mathbb R^+$, write
$$u(x)=u_+(x)+u_{-}(x)$$
such that $u_+(x)$ is supported on $[1, \infty)$ and $u_{-}(x)$ is supported on $(0,1]$. 
Let $f(x)$ be a function in $L^2(\mathbb R^+, (1+x^{\epsilon}) d x)$. Then Theorem \ref{fund} asserts that
$$ \int_0^{T_1} |\mu_{\f d, \exp(-2\pi i x T) f}( a)_+|^2_{L^2(\mathbb R^+, a^{\epsilon} \frac{d a}{a})} d T  \leq 2 T_1 C_{\epsilon} \| f_+ \|^2_{L^2(\mathbb R^+, a^{\epsilon} \frac{d a}{a})} .$$
By Fubini's theorem, for almost all $T \in [0,T_1]$, $\mu_{\f d, \exp(-2\pi i x T) f}( a)_+ \in L^2(\mathbb R^+, a^{\epsilon} \frac{d a}{a})$. \\
\\
Recall that for $0 < \sigma \leq \frac{\epsilon}{2}$, The Fourier-Mellin transform $\mc M^{\sigma} f$ is well-defined and
$$\| \mc M^{\sigma} f \|^2= \frac{1}{2 \pi} \int_{\mathbb R^+} |f(x)|^2 x^{2 \sigma-1} d x$$
We can now restate Theorem \ref{fund}.
\begin{thm} Let $L(s, \f d) \in \mc R$,  $T_1 \geq 1$ and $\epsilon >0$. Let $f \in L^2(\mathbb R^+, (1+x^{\epsilon}) d x)$. There exists a constant $C_{\epsilon, T_1, \f d}$ such that
$$ \int_0^{T_1} \|\mc M^{\frac{\epsilon}{2}} (\mu_{\f d, \exp (-2 \pi i T x) f}( a)_+)(s_1)\|_{L^2(\mathbb R)}^2 d T  \leq  C_{\epsilon, T_1, \f d} \| \mc M^{\frac{1+\epsilon}{2}}(f_+)(s_1) \|_{L^2(\mathbb R)}^2.$$
\end{thm}
Recall that $$\mc M^{\frac{\epsilon}{2}} (\mu_{\f d, \exp (-2 \pi i T x) f}( a))(s_1)=L(\frac{1+\epsilon}{2}+ s_1 i, \f d) \mc M^{\frac{1+\epsilon}{2}} \{f(x) \exp - 2\pi i T x \} (s_1)$$
The Mellin transform $\mc M \{\mu_{\f d, \exp (-2 \pi i T x) f}( a)_+ \}(\frac{\epsilon}{2}+s_1 i)$, which lives in the Hardy space of a right half plane, may be obtained explicitly from The Mellin transform
$$\mc M \{\mu_{\f d, \exp (-2 \pi i T x) f}( a)\}(\frac{\epsilon}{2}+s_1 i)$$
 by solving a Riemann-Hilbert problem. We indeed get an inequlity in terms of the Dirichlet series $L(\frac{1+\epsilon}{2}+ s_1 i, \f d)$. Nevertheless, in order to obtain a bound on $L(s, \f d)$ we will have to control 
$$ \int_0^{T_1} \|\mc M^{\frac{\epsilon}{2}} (\mu_{\f d, \exp (-2 \pi i T x) f}( a)_-)(s_1)\|_{L^2(\mathbb R)}^2 d T $$
with $ a \in (0,1]$. This can be achieved when $\mathfrak d$ comes from an automorphic distribution for an arithmetic subgroup of $SL(2, \mathbb Z)$ (\cite{sc} \cite{be} \cite{he1}).


\begin{thebibliography}{99}

\bibitem{be}  J. Bernstein and A. Reznikov, Sobolev norms of automorphic functionals and Fourier coefficients of cusp forms. C. R. Acad. Sci. Paris Ser. I Math. 327 (1998), no. 2, 111-116. 
\bibitem{borel} A. Borel {\it Automorphic forms on $SL(2)$}  Cambridge Tracts in Mathematics, 130. Cambridge University Press, Cambridge, 1997. 

\bibitem{fo} G. Folland {\it A course in abstract harmonic analysis.}
Studies in Advanced Mathematics. CRC Press, Boca Raton, FL, 1995.
\bibitem{he}  H. He Generalized matrix coefficients for infinite dimensional unitary representations. J. Ramanujan Math. Soc. 29 (2014), no. 3, 253-272.
\bibitem{he1} H. He Certain $L^2$-norm on Automorphic Representations  of $SL(2, \mathbb R)$ (preprint).

\bibitem{knapp} A. Knapp {\it Representation theory of semisimple Groups} Princeton University Press 2002.
\bibitem{sc} W. Schmid, Automorphic distributions for SL(2,R). Conference Moshe Flato 1999, Vol. I (Dijon), 345-387, Math. Phys. Stud., 21, Kluwer Acad. Publ., Dordrecht, 2000. 
\bibitem{ti} E. Titchmarsh, Introduction to the theory of Fourier Integrals, Oxford at the Clarendon Press (1962).
\bibitem{tr}  F. Treves, Topological Vector Spaces, Distributions and
Kernels, Academic Press, New York 1967.
\bibitem{wa} N. Wallach {\it Real Reductive Groups I, II} Academic Press 1992.
\bibitem{war}  G. Warner, {\it Harmonic Analysis on Semi-Simple Lie
Groups I}, Springer-Verlag, 1972
\bibitem{zyg} A. Zygmund,  {\it Trigonometric Series}, Volume 2, Cambridge University Press 1959.
\end{thebibliography}
\end{document}